\documentclass[a4paper,12pt]{article}
\usepackage[all]{xy}
\usepackage[T1]{fontenc}
\usepackage[dvips]{graphicx}
\usepackage{amsmath,amssymb,amsfonts,amsthm,latexsym,mathrsfs,textcomp,verbatim}
\xyoption{web}

\title{Fra\"iss\'e's construction\\ from a topos-theoretic perspective}
\author{{Olivia Caramello} \vspace{3 mm}\\ {\small DPMMS, University of Cambridge,}\\{\small Wilberforce Road, Cambridge CB3 0WB, UK}\\{\small O.Caramello@dpmms.cam.ac.uk}}
\date{\today}
\begin{document}
\bgroup           % fake a titlepage 
\let\footnoterule\relax  % no rule above thanks footnotes 
\maketitle
\flushleft  
\begin{abstract}
We present a topos-theoretic interpretation of (a categorical generalization of) Fra\"iss\'e's construction in model theory, with applications to countably categorical theories. 
\end{abstract} 
\egroup 
\flushleft
\vspace{5 mm}

%MACROS-----------------------------------------------------------------------------------------------------------------------

%	European dates ``19 April 1990'' not ``April 19, 1990''
\def\Monthnameof#1{\ifcase#1\or
   January\or February\or March\or April\or May\or June\or
   July\or August\or September\or October\or November\or December\fi}
\def\today{\number\day~\Monthnameof\month~\number\year}

%===========================================================================
%	END OF PROOF BOX
%
%
%  The complexity of the macro necessary to get a little box on the
%  right-hand-side at the end of a proof is amazing.  It really does
%  have to be this long!  Otherwise you're liable to get it at the
%  beginning of the next line, or even on the next page.
%
\def\pushright#1{{%        set up
   \parfillskip=0pt            % so \par doesnt push \square to left
   \widowpenalty=10000         % so we dont break the page before \square
   \displaywidowpenalty=10000  % ditto
   \finalhyphendemerits=0      % TeXbook exercise 14.32
  %
  %                 horizontal
   \leavevmode                 % \nobreak means lines not pages
   \unskip                     % remove previous space or glue
   \nobreak                    % don't break lines
   \hfil                       % ragged right if we spill over
   \penalty50                  % discouragement to do so
   \hskip.2em                  % ensure some space
   \null                       % anchor following \hfill
   \hfill                      % push \square to right
   {#1}                        % the end-of-proof mark (or whatever)
  %
  %                   vertical
   \par}}                      % build paragraph

% prefer proofs with statements, also space after
\def\qed{\pushright{$\square$}\penalty-700 \smallskip}

\newtheorem{theorem}{Theorem}[section]

\newtheorem{proposition}[theorem]{Proposition}

\newtheorem{scholium}[theorem]{Scholium}

\newtheorem{lemma}[theorem]{Lemma}

\newtheorem{corollary}[theorem]{Corollary}

\newtheorem{conjecture}[theorem]{Conjecture}

\newenvironment{proofs}%
 {\begin{trivlist}\item[]{\bf Proof }}%
 {\qed\end{trivlist}}

  \newtheorem{rmk}[theorem]{Remark}
\newenvironment{remark}{\begin{rmk}\em}{\end{rmk}}

  \newtheorem{rmks}[theorem]{Remarks}
\newenvironment{remarks}{\begin{rmks}\em}{\end{rmks}}

  \newtheorem{defn}[theorem]{Definition}
\newenvironment{definition}{\begin{defn}\em}{\end{defn}}

  \newtheorem{eg}[theorem]{Example}
\newenvironment{example}{\begin{eg}\em}{\end{eg}}

  \newtheorem{egs}[theorem]{Examples}
\newenvironment{examples}{\begin{egs}\em}{\end{egs}}

%%%%%%%%%%%%%%%%%%%%%%%%%%%%%%%%%%%%%%%%%%%%%%%%%%%%%%%%%%%%%%%%%%%%%%

%  change some single-character symbols to be more appropriate for logic
\mathcode`\<="4268  % < = \langle
\mathcode`\>="5269  % > = \rangle
\mathcode`\.="313A  % make . a binary  relation
\mathchardef\semicolon="603B % the original
\mathchardef\gt="313E
\mathchardef\lt="313C

\newcommand{\app}% application
 {{\sf app}}

\newcommand{\Ass}% category of assemblies
 {{\bf Ass}}

\newcommand{\ASS}% indexed version of \Ass
 {{\mathbb A}{\sf ss}}

\newcommand{\Bb}%blackboard bold
{\mathbb}

\newcommand{\biimp}% bi-implication
 {\!\Leftrightarrow\!}

\newcommand{\bim}% bimorphism
 {\rightarrowtail\kern-1em\twoheadrightarrow}

\newcommand{\bjg}% bi-judgement
 {\mathrel{{\dashv}\,{\vdash}}}

\newcommand{\bstp}[3]% bimorphism version of \stp
 {\mbox{$#1\! : #2 \bim #3$}}

\newcommand{\cat}%concatenation
 {\!\mbox{\t{\ }}}

\newcommand{\cinf}%C-infinity
 {C^{\infty}}

\newcommand{\cinfrg}%category of C-infinity rings
 {\cinf\hy{\bf Rng}}

\newcommand{\cocomma}[2]% cocomma object
 {\mbox{$(#1\!\uparrow\!#2)$}}

\newcommand{\cod}% codomain
 {{\rm cod}}

\newcommand{\comma}[2]% comma object
 {\mbox{$(#1\!\downarrow\!#2)$}}

\newcommand{\comp}% composition
 {\circ}

\newcommand{\cons}% concatenation
 {{\sf cons}}

\newcommand{\Cont}% category of continuous G-sets
 {{\bf Cont}}

\newcommand{\ContE}% continuous G-sets relative to $\cal E$
 {{\bf Cont}_{\cal E}}

\newcommand{\ContS}% ditto for $\cal S$
 {{\bf Cont}_{\cal S}}

\newcommand{\cover}% cover
 {-\!\!\triangleright\,}

\newcommand{\cstp}[3]% cover version of \stp
 {\mbox{$#1\! : #2 \cover #3$}}

\newcommand{\Dec}% decalage
 {{\rm Dec}}

\newcommand{\DEC}% decalage (\Bbb version)
 {{\mathbb D}{\sf ec}}

\newcommand{\den}[1]% denotation of #1 
 {[\![#1]\!]}

\newcommand{\Desc}% category of descent data
 {{\bf Desc}}

\newcommand{\dom}% domain
 {{\rm dom}}

\newcommand{\Eff}% effective topos
 {{\bf Eff}}

\newcommand{\EFF}% indexed version of \Eff
 {{\mathbb E}{\sf ff}}

\newcommand{\empstg}% empty string
 {[\,]}

\newcommand{\epi}% epimorphism
 {\twoheadrightarrow}

\newcommand{\estp}[3]% epimorphism version of \stp
 {\mbox{$#1 \! : #2 \epi #3$}}

\newcommand{\ev}% evaluation
 {{\rm ev}}

\newcommand{\Ext}% category of extracts
 {{\rm Ext}}

\newcommand{\fr}% Fraktur (i.e. Gothic)
 {\sf}

\newcommand{\fst}% first projection
 {{\sf fst}}

\newcommand{\fun}[2]% function-type
 {\mbox{$[#1\!\to\!#2]$}}

\newcommand{\funs}[2]% function-type as subscript
 {[#1\!\to\!#2]}

\newcommand{\Gl}% topos obtained by glueing
 {{\bf Gl}}

\newcommand{\hash}% hash sign (used as infix)
 {\,\#\,}

\newcommand{\hy}% hyphen (in math mode)
 {\mbox{-}}

\newcommand{\im}% image
 {{\rm im}}

\newcommand{\imp}% implication
 {\!\Rightarrow\!}

\newcommand{\Ind}[1]% ind-completion of #1
 {{\rm Ind}\hy #1}

\newcommand{\iten}[1]% enumerated item
{\item[{\rm (#1)}]}

\newcommand{\iter}% iterator
 {{\sf iter}}

\newcommand{\Kalg}%category of $K$-algebras
 {K\hy{\bf Alg}}

\newcommand{\llim}% left (inverse) limit
 {{\mbox{$\lower.95ex\hbox{{\rm lim}}$}\atop{\scriptstyle
{\leftarrow}}}{}}

\newcommand{\llimd}% \displaymath version of \llim
 {\lower0.37ex\hbox{$\pile{\lim \\ {\scriptstyle
\leftarrow}}$}{}}

\newcommand{\Mf}%category of manifolds
 {{\bf Mf}}

\newcommand{\Mod}% category of modest assemblies
 {{\bf Mod}}

\newcommand{\MOD}% indexed version of \Mod
{{\mathbb M}{\sf od}}

\newcommand{\mono}% monomorphism 
 {\rightarrowtail}

\newcommand{\mor}% class of morphisms
 {{\rm mor}}

\newcommand{\mstp}[3]% monomorphism version of \stp
 {\mbox{$#1\! : #2 \mono #3$}}

\newcommand{\Mu}%capital mu
 {{\rm M}}

\newcommand{\name}[1]% name of a relation
 {\mbox{$\ulcorner #1 \urcorner$}}

\newcommand{\names}[1]%\name used as subscript
 {\mbox{$\ulcorner$} #1 \mbox{$\urcorner$}}

\newcommand{\nml}% normal subgroup
 {\triangleleft}

\newcommand{\ob}% class of objects
 {{\rm ob}}

\newcommand{\op}% opposite category
 {^{\rm op}}

\newcommand{\pepi}% partial epimorphism
 {\rightharpoondown\kern-0.9em\rightharpoondown}

\newcommand{\pmap}% partial map arrow
 {\rightharpoondown}

\newcommand{\Pos}% positivization of a coherent category
 {{\bf Pos}}

\newcommand{\prarr}% parallel pair of arrows
 {\rightrightarrows}

\newcommand{\princfil}[1]% principal filter
 {\mbox{$\uparrow\!(#1)$}}

\newcommand{\princid}[1]% principal ideal
 {\mbox{$\downarrow\!(#1)$}}

\newcommand{\prstp}[3]% parallel-pair version of \stp
 {\mbox{$#1\! : #2 \prarr #3$}}

\newcommand{\pstp}[3]% partial-map version of \stp
 {\mbox{$#1\! : #2 \pmap #3$}}

\newcommand{\relarr}% relation-type arrow
 {\looparrowright}

\newcommand{\rlim}% right limit, i.e. colimit
 {{\mbox{$\lower.95ex\hbox{{\rm lim}}$}\atop{\scriptstyle
{\rightarrow}}}{}}

\newcommand{\rlimd}% \displaymath version of \rlim
 {\lower0.37ex\hbox{$\pile{\lim \\ {\scriptstyle
\rightarrow}}$}{}}

\newcommand{\rstp}[3]% relation version of \stp
 {\mbox{$#1\! : #2 \relarr #3$}}

\newcommand{\scn}% Sierpinski cone
 {{\bf scn}}

\newcommand{\scnS}% ditto relative to $\cal S$
 {{\bf scn}_{\cal S}}

\newcommand{\semid}% semidirect product
 {\rtimes}

\newcommand{\Sep}% category of separated objects
 {{\bf Sep}}

\newcommand{\sep}% category of separated objects
 {{\bf sep}}

\newcommand{\Set}% category of sets
 {{\bf Set }}

\newcommand{\Sh}% category of sheaves
 {{\bf Sh}}

\newcommand{\ShE}% sheaves relative to $\cal E$
 {{\bf Sh}_{\cal E}}

\newcommand{\ShS}% ditto for $\cal S$
 {{\bf Sh}_{\cal S}}

\newcommand{\sh}% category of sheaves
 {{\bf sh}}

\newcommand{\Simp}% the simplicial category
 {{\bf \Delta}}

\newcommand{\snd}% second projection
 {{\sf snd}}

\newcommand{\stg}[1]% string of #1
 {\vec{#1}}

\newcommand{\stp}[3]% source-target predicate
 {\mbox{$#1\! : #2 \to #3$}}

\newcommand{\Sub}% subobject lattice
 {{\rm Sub}}

\newcommand{\SUB}% indexed category of subobjects
 {{\mathbb S}{\sf ub}}

\newcommand{\tbel}% totally below
 {\prec\!\prec}

\newcommand{\tic}[2]%term-in-context, etc.
 {\mbox{$#1\!.\!#2$}}

\newcommand{\tp}% is of type
 {\!:}

\newcommand{\tps}% subscript version of \tp
 {:}

\newcommand{\tsub}% truncated subtraction
 {\pile{\lower0.5ex\hbox{.} \\ -}}

\newcommand{\wavy}% wavy arrow
 {\leadsto}

\newcommand{\wavydown}% wavy downarrow
 {\,{\mbox{\raise.2ex\hbox{\hbox{$\wr$}
\kern-.73em{\lower.5ex\hbox{$\scriptstyle{\vee}$}}}}}\,}

\newcommand{\wbel}% way-below relation
 {\lt\!\lt}

\newcommand{\wstp}[3]% wavy version of \stp
 {\mbox{$#1\!: #2 \wavy #3$}}

\newcommand{\fu}[2]
{[#1,#2]}

\newcommand{\st}[2]% source-target predicate
 {\mbox{$#1 \to #2$}}

\flushleft
\section{A categorical generalization of Fra\"iss\'e's theorem}
In this section we present a categorical generalization of Fra\"iss\'e's construction in model theory.
Our result is technically similar to (though more general than) the categorical theorem in \cite{DG1}, but follows as an application of the theory developed by Kubi\'s in \cite{Kubis}.\\
First, let us introduce the relevant terminology.
\begin{definition}
A category $\cal C$ is said to satisfy the \emph{amalgamation property} (AP) if for every objects $a,b,c\in {\cal C}$ and morphisms $f:a\rightarrow b$, $g:a\rightarrow c$ in $\cal C$ there exists an object $d\in \cal C$ and morphisms $f':b\rightarrow d$, $g':c\rightarrow d$ in $\cal C$ such that $f'\circ f=g'\circ g$:
\[  
\xymatrix {
a \ar[d]_{g} \ar[r]^{f} & b  \ar@{-->}[d]^{f'} \\
c \ar@{-->}[r]_{g'} & d } 
\] 
\end{definition} 
Notice that $\cal C$ satisfies the amalgamation property if and only if ${\cal C}^{\textrm{op}}$ satisfies the right Ore condition. So if $\cal C$ satisfies AP then we may equip ${\cal C}^{\textrm{op}}$ with the atomic topology. This point will be the basis of our topos-theoretic interpretation described in the next section. 
\begin{definition}
A category $\cal C$ is said to satisfy the \emph{joint embedding property} (JEP) if for every pair of objects $a,b\in {\cal C}$ there exists an object $c\in \cal C$ and morphisms $f:a\rightarrow c$, $g:b\rightarrow c$ in $\cal C$:
\[  
\xymatrix {
 & a  \ar@{-->}[d]^{f} \\
b \ar@{-->}[r]_{g} & c } 
\] 
\end{definition} 
Notice that if $\cal C$ has a weakly initial object then AP on $\cal C$ implies JEP on $\cal C$; however in general the two notions are quite distinct from one another. 
\begin{definition}
Given an embedding $i:{\cal C}\rightarrow {\cal D}$, an object $u\in {\cal D}$ is said to be \emph{$\cal C$-homogeneous} if for every objects $a,b \in {\cal C}$ and arrows $j:a\rightarrow b$ in ${\cal C}$ and $\chi:a\rightarrow u$ in $\cal D$ there exists an arrow $\tilde{\chi}:b\rightarrow u$ such that $\tilde{\chi}\circ j=\chi$:
\[  
\xymatrix {
a \ar[d]_{j} \ar[r]^{\chi} & u \\
b \ar@{-->}[ur]_{\tilde{\chi}} &  } 
\] 
$u$ is said to be \emph{$\cal C$-ultrahomogeneous} if for every objects $a,b \in {\cal C}$ and arrows $j:a\rightarrow b$ in ${\cal C}$ and $\chi_{1}:a\rightarrow u$, $\chi_{2}:b\rightarrow u$ in $\cal D$ there exists an isomorphism $\check{j}:u\rightarrow u$ such that $\check{j}\circ \chi_{1}=\chi_{2}\circ j$:
\[  
\xymatrix {
a \ar[d]_{j} \ar[r]^{\chi_{1}} & u \ar@{-->}[d]^{\check{j}}\\
b \ar[r]_{\chi_{2}} & u } 
\] 
$u$ is said to be \emph{$\cal C$-universal} if it is $\cal C$-cofinal, that is for every $a\in {\cal C}$ there exists an arrow $\chi:a\rightarrow u$ in $\cal D$:
\[  
\xymatrix {
a \ar@{-->}[r]^{\chi} & u  } 
\]     
\end{definition} 
\begin{rmks}
\emph{It is easy to see that if $u$ is $\cal C$-ultrahomogeneous and $\cal C$-universal then $u$ is $\cal C$-homogeneous. Also, to verify that an object $u$ in $\cal D$ is $\cal C$-ultrahomogeneous one can clearly suppose, without loss of generality, that the arrow $j$ in the definition is an identity.} 
\end{rmks}
Let us recall the following definitions from \cite{Kubis}.\\
Given a category $\cal C$ and a collection of arrows ${\cal F}\subseteq arr(\cal C)$, $\cal F$ is said to be dominating in $\cal C$ if the family $Dom(\cal F)$ of objects which are domains of an arrow in $\cal F$ is cofinal in $\cal C$ and satisfies the following property: for every $a\in Dom(\cal F)$ and every arrow $f:a\rightarrow x$ in $\cal C$ there exists an arrow $g:x\rightarrow cod(g)$ in $\cal C$ such that $g\circ f\in \cal F$.\\
Notice that $arr(\cal C)$ is always dominating in $\cal C$, and if $\cal C'$ is a skeleton of $\cal C$, $arr(\cal C')$ is dominating in $\cal C$.\\
Given a category $\cal C$ and an ordinal $\kappa\gt 0$, an inductive $\kappa$-sequence (or $\kappa$-chain) in $\cal C$ is a functor $\vec{u}:\kappa\rightarrow \cal C$, where $\kappa$ is regarded as a poset category. For $i\in \kappa$ we denote $\vec{u}(i)$ by $u_{i}$ and for $i,j\in \kappa$ such that $i\leq j$ we denote $\vec{u}(i\rightarrow j):u_{i}\rightarrow u_{j}$ by $u_{i}^{j}$. $\vec{u}$ is said to be a Fra\"iss\'e sequence of length $\kappa$ (or, briefly, a $\kappa$-Fra\"iss\'e sequence) in $\cal C$ if it satisfies the following conditions:\\
(1) For every $a\in \cal C$ there exists $i\in \kappa$ and an arrow $\chi:a\rightarrow u_{i}$ in $\cal C$;\\
(2) For every $i\in \kappa$ and for every arrow $f:u_{i}\rightarrow cod(f)$ in $\cal C$, there exists $j \in \kappa$ with $j\geq i$ and an arrow $g:cod(f)\rightarrow u_{j}$ such that $u_{i}^{j}=g\circ f$.\\  
$\vec{u}$ is said to have the extension property if it satisfies the following condition:\\
For every arrows $f:a\rightarrow b, g:a\rightarrow u_{i}$ in $\cal C$ where $i\in \kappa$, there exists $j \in \kappa$ with $j\geq i$ and an arrow $h:b\rightarrow u_{j}$ such that $u_{i}^{j}\circ g = h\circ f$.\\
Of course, every sequence satisfying the extension property satisfies property (2) in the definition of Fra\"iss\'e sequence.\\ 
A category $\cal C$ is said to be $\kappa$-bounded if every chain in $\cal C$ of length $\lambda\lt \kappa$ has a cocone in $\cal C$ over it. Clearly, every category is $\omega$-bounded.\\  
A $\kappa$-chain $\vec{u}:\kappa\rightarrow \cal C$ is said to be continuous if for each limit ordinal $j\in \kappa$, $u_{j}$ is the colimit of the $j$-chain obtained as the restriction of $\vec{u}$ to $j$ with universal colimit arrows given by the arrows $u_{i}\rightarrow u_{j}$ ($i\lt j$) of the chain.\\   
Given an infinite cardinal $\kappa$ and an embedding $i:{\cal C}\rightarrow {\cal D}$, we denote by ${\cal D}_{\kappa}$ the full subcategory of $\cal D$ on the objects that can be expressed as colimits of $\kappa$-chains in $\cal C$ and by ${\cal D}_{\kappa}^{c}$ the full subcategory of $\cal D$ on the objects that can be expressed as colimits of continuous $\kappa$-chains in $\cal C$. We will say that an embedding $i:{\cal C}\rightarrow {\cal D}$ is $\kappa$-continuous if ${\cal D}_{\kappa}={\cal D}_{\kappa}^{c}$. Obviously, every embedding is $\omega$-continuous.\\ 
Following the terminology in \cite{DG1}, we will say that an object $a$ in $\cal C$ is $\kappa$-small in $\cal D$ if the functor $Hom_{\cal D}(i(a),- ):{\cal D}\rightarrow \Set$ preserves all colimits of $\kappa$-chains in $\cal D$; in particular, every finitely presentable object in $\cal C$ is  $\kappa$-small.\\
Notice that, given an embedding $i:{\cal C}\rightarrow {\cal D}$ such that all the objects in $\cal C$ are $\kappa$-small in $\cal D$, for $i$ to be $\kappa$-continuous it suffices that $\cal C$ is closed under colimits of $\lambda$-chains in $\cal D$ for each $\lambda \lt \kappa$; indeed, given an inductive $\kappa$-sequence $\vec{u}$ in $\cal C$ with colimit $u$ we can construct (by transfinite recursion) a continuous $\kappa$-chain $\vec{v}$ in $\cal C$ with a universal colimiting cone $D$ to $u$ (cfr. also the proof of Lemma 1 in \cite{rosicky}); indeed, denoted by $j_{i}:u_{i}\rightarrow u$ (for $i\lt \kappa$) the universal colimit arrows for $\vec{u}$, we define $\vec{v}$ as follows:\\
$\vec{v}(0)=\vec{u}(0)$ and $D(0)=j_{0}$;\\   
given $\vec{v}(i)$ and $D(i):v_{i}\rightarrow u$, $v_{i}$ being $\kappa$-small in $\cal D$, there exists $j\gt i$ and an arrow $h:v_{i}\rightarrow u_{j}$ such that $D(i)=j_{i}\circ h$; we put $\vec{v}(i+1)=u_{j}$ and $D(i+1)=h$;\\
if $i\lt \kappa$ is a limit ordinal then we define $\vec{v}(i)$ and $D(i)$ respectively as the colimit $colim_{j\lt i}\vec{v}(j)$ and the unique arrow $colim_{j\lt i}\vec{v}(j)\rightarrow u$ induced via the universal property of the colimit by the arrows $D(j):\vec{v}(j)\rightarrow u$ (for $j\lt i$).\\
The sequence $\vec{v}$ is defined on the arrows in the obvious way.\\
If $i$ is the embedding of the full subcategory on the $\kappa$-presentable objects of a $\kappa$-accessible category $\cal C$ having directed colimits into $\cal C$, then, denoted by ${\cal C}^{{\kappa}^{+}}$ the full subcategory of $\cal C$ on the ${\kappa}^{+}$-presentable objects, by the proof of Lemma 1 in \cite{rosicky} we have that ${\cal C}^{{\kappa}^{+}}={\cal C}_{\kappa}={\cal C}_{\kappa}^{c}$; in particular $i$ is $\kappa$-continuous.   
\begin{theorem}\label{teofond}
Let $\kappa$ be an infinite regular cardinal and $\cal C$ be a $\kappa$-bounded category satisfying the amalgamation and the joint embedding properties. If there exists a dominating family of arrows $\cal F$ in $\cal C$ such that $|{\cal F}|\leq \kappa$, then for any embedding $i:{\cal C}\rightarrow D$ such that $\cal D$ has all colimits of $\kappa$-chains in $\cal C$ and all the objects in $\cal C$ are $\kappa$-small in $\cal D$, there exists in ${\cal D}_{\kappa}$ a $\cal C$-homogeneous and $\cal C$-universal object; if moreover all the morphisms in ${\cal D}_{\kappa}^{c}$ are monic (as arrows in ${\cal D}_{\kappa}^{c}$) then every $\cal C$-homogeneous and $\cal C$-universal object in ${\cal D}_{\kappa}^{c}$ is $\cal C$-ultrahomogeneous and unique (up to isomorphism) with these properties in ${\cal D}_{\kappa}^{c}$.\\
Conversely, given an embedding $i:{\cal C}\rightarrow {\cal D}$ such that all the morphisms in ${\cal D}_{\kappa}$ are monic, if there exists in ${\cal D}_{\kappa}$ an object which is $\cal C$-homogeneous and $\cal C$-universal, then the category $\cal C$ satisfies the amalgamation and joint embedding properties.      
\end{theorem}
\begin{proofs}
Let $u$ be the colimit in $\cal D$ of an inductive $\kappa$-sequence $\vec{u}$ in $\cal C$. Then the following facts hold:\\
(1) If $u$ is $\cal C$-homogeneous and $\cal C$-universal and all the morphisms in ${\cal D}_{\kappa}$ (respectively, in ${\cal D}_{\kappa}^{c}$ if $u$ belongs to ${\cal D}_{\kappa}^{c}$) are monic then $\vec{u}$ is a Fra\"iss\'e sequence.\\
(2) If $\vec{u}$ is a $\kappa$-Fra\"iss\'e sequence then $u$ is $\cal C$-homogeneous and $\cal C$-universal; moreover, if $\vec{u}$ is continuous then $u$ is $\cal C$-ultrahomogeneous.\\
To prove (1), let us suppose that $u$ is $\cal C$-homogeneous and $\cal C$-universal and all the morphisms in ${\cal D}_{\kappa}$ (respectively, in ${\cal D}_{\kappa}^{c}$ if $u$ belongs to ${\cal D}_{\kappa}^{c}$) are monic. Condition (1) in the definition of Fra\"iss\'e sequence trivially follows from the fact that $u$ is $\cal C$-universal and every object of $\cal C$ is $\kappa$-small in ${\cal D}$. To verify condition (2), we prove that $\vec{u}$ satisfies the extension property.  Since $\vec{u}$ is $\cal C$-homogeneous, then given arrows $f:a\rightarrow b$ and $g:a\rightarrow u_{i}$ in $\cal C$ where $i\in \kappa$, and the colimit map $j_{i}:u_{i}\rightarrow u$, there exists an arrow $h:b\rightarrow u$ such that $h\circ f=j_{i}\circ g$. Now, $b$ being $\kappa$-small in ${\cal D}$, $h$ factors as $b\stackrel{h_{j}}{\rightarrow}u_{j}\stackrel{j_{j}}{\rightarrow}u$ for a sufficiently large $j$. If we take a $j\geq i$ then we clearly have $u_{i}^{j}\circ g=h_{j}\circ f$, $j_{j}$ being monic.\\
Let us now prove fact (2). By condition (1) in the definition of Fra\"iss\'e sequence, $u$ is clearly $\cal C$-universal. Let us now prove that $u$ is $\cal C$-homogeneous.\\
Given objects $a,b \in {\cal C}$ and arrows $f:a\rightarrow b$ in ${\cal C}$ and $\chi:a\rightarrow u$ in $\cal D$, we want to prove that there exists an arrow $\tilde{\chi}:b\rightarrow u$ such that $\tilde{\chi}\circ f=\chi$:
\[  
\xymatrix {
a \ar[d]_{f} \ar[r]^{\chi} & u \\
b \ar@{-->}[ur]_{\tilde{\chi}} &  } 
\] 
Since $a$ is $\kappa$-small in $\cal D$, $\chi$ factors as $a\stackrel{\chi_{i}}{\rightarrow}u_{i}\stackrel{j_{i}}{\rightarrow}u$ for some $i\in \kappa$. From the fact that $\cal C$ satisfies AP we obtain an object $d\in {\cal C}$ and two arrows $h:u_{i}\rightarrow d$ and $l:b\rightarrow d$ such that $l\circ f=h\circ \chi_{i}$. Now by condition (2) in the definition of Fra\"iss\'e sequence we get a $j\in \kappa$ with $j\geq i$ and an arrow $m:d\rightarrow u_{j}$ such that $u_{i}^{j}=m\circ h$. Hence the arrow $\tilde{\chi}:=j_{j}\circ m\circ l$ satisfies the required property.\\  
Let us now prove the following fact, to which we will refer as to fact (3):\\
If $\vec{u}$ and $\vec{v}$ are two continuous $\kappa$-Fra\"iss\'e sequences in $\cal C$ and $f:u_{k}\rightarrow v_{l}$ is an arrow between ``elements'' respectively of $\vec{u}$ and $\vec{v}$ there exists in $\cal D$ an isomorphism $\tilde{f}:\varinjlim \vec{u}\rightarrow \varinjlim \vec{v}$ such that $\tilde{f}\circ j_{k}=j'_{l}\circ f$ (where $j_{k}:u_{k}\rightarrow \varinjlim \vec{u}$ and $j'_{l}:v_{l}\rightarrow \varinjlim \vec{v}$ are the obvious colimit arrows).\\
To this end, let us establish the following fact: given any $k, l \in \kappa$ and any arrow $f:u_{k}\rightarrow v_{l}$ there exist two strictly increasing functions $k,l:\kappa\rightarrow \kappa$ and two natural transformations $F:\vec{u}\circ k\rightarrow \vec{v}\circ l$ and $G:\vec{v}\circ l\rightarrow \vec{u}\circ k^{+}$, where $k^{+}$ is the function defined by $k^{+}(i)=k(i+1)$ for each $i\lt \kappa$, with the following properties:\\
$k(0)=k$ and $l(0)=l$,\\ 
$F(0)=f$ and $F(i+1)\circ G(i)=v^{l(i+1)}_{l(i)}, G(i)\circ F(i)=u^{k(i+1)}_{k(i)}$ (for each $i\in \kappa$).\\
We define $k(i), l(i), F(i), G(i)$ (and prove that they satisfy the required properties) by transfinite induction on $i\lt \kappa$.\\
For $i=0$ we put $k(0)=k$, $l(0)=l$, $F(0)=f$ and define $k(1)$ and $G(0):v_{l}\rightarrow u_{k(1)}$ as follows: by condition (2) in the definition of Fra\"iss\'e sequence applied to $\vec{u}$ there exist $j\in \kappa$ with $j\gt k$ and an arrow $s:v_{l}\rightarrow u_{j}$ such that $s\circ f=u_{k}^{j}$; we put $k(1)=j$ and $G(0)=s$.\\     
Given $k(i), k(i+1), l(i), F(i), G(i)$ we define $k(i+2), l(i+1), F(i+1), G(i+1)$ as follows: by condition (2) in the definition of Fra\"iss\'e sequence applied to $\vec{v}$ there exist $j\in \kappa$ with $j\gt l(i)$ and an arrow $s:u_{k(i+1)}\rightarrow v_{j}$ such that $s\circ G(i)=v_{l(i)}^{j}$; we put $l(i+1)=j$ and $F(i+1)=s$. Again, by condition (2) in the definition of Fra\"iss\'e sequence applied to $\vec{u}$ there exist $j'\in \kappa$ with $j'\gt k(i+1)$ and an arrow $s':v_{l(i+1)}\rightarrow u_{j'}$ such that $s'\circ F(i+i)=u_{k(i+1)}^{j'}$; we put $k(i+2)=j'$ and $G(i+1)=s'$.\\
If $i=sup_{j\lt i}j$ is a limit ordinal we put $k(i)=sup_{j\lt i}k(j)\in \kappa$, $l(i)=sup_{j\lt i}l(j)\in \kappa$ and define both $F(i)$ and $G(i)$ by taking colimits. More precisely, since the restriction of $k$ to $i$ is strictly increasing and the chain $\vec{u}$ is $\kappa$-continuous then $u_{k(i)}$ is the colimit of the restriction of the chain $\vec{u}\circ k$ to $i$; analogously, $v_{l(i)}$ is the colimit of the restriction of the chain $\vec{v}\circ l$ to $i$; then we define $F(i):u_{k(i)}\rightarrow v_{l(i)}$ to be the unique arrow, given by the universal property of the colimit, such that for each $j\lt i$ $F(i)\circ u_{k(j)}^{k(i)}=v_{l(j)}^{l(i)}\circ F(j)$. $G(i)$ is defined similarly.\\ The verification that all the required properties are satisfied is easily done by induction on $i\lt k$.\\    
Now, since $k$ and $l$ are strictly increasing functions, then, regarded as functors $\kappa \rightarrow \kappa$, they are cofinal. This implies that $u=\varinjlim \vec{u}=\varinjlim (\vec{u}\circ k)$ and $v=\varinjlim \vec{v}=\varinjlim (\vec{v}\circ l)$. Hence the natural transformations $F$ and $G$ respectively induce arrows $\tilde{f}:u\rightarrow v$ and $g:v\rightarrow u$ such that, denoted by $j_{k(i)}:u_{k(i)}\rightarrow u$ and $j'_{l(i)}:v_{l(i)}\rightarrow v$ the colimit arrows, $\tilde{f}\circ j_{k(i)}=j'_{l(i)}\circ F(i)$ and $g\circ j'_{l(i)}=j_{k(i+1)}\circ G(i)$ for each $i\in \kappa$; in particular $\tilde{f}\circ j_{k}=j'_{l}\circ f$. We have $g\circ \tilde{f}=1_{u}$ and $\tilde{f}\circ g=1_{v}$ in $\cal D$, from which it follows that $\tilde{f}$ is an isomorphism with the required property. Indeed, let us for example prove that first equality; the second follows similarly. By the universal property of the colimit $u=\varinjlim (\vec{u}\circ k)$, it is equivalent to check that $g\circ f\circ j_{k(i)}=j_{k(i)}$ for each $i\in \kappa$. Now by the equalities above we obtain $g\circ \tilde{f}\circ j_{k(i)}=g\circ j'_{l(i)}\circ F(i)=j_{k(i+1)}\circ G(i)\circ F(i)=j_{k(i+1)}\circ u^{k(i+1)}_{k(i)}=j_{k(i)}$, as required. This completes the proof of fact (3).\\
Coming back to our Fra\"iss\'e sequence $\vec{u}$, by taking $\vec{v}=\vec{u}$ in fact (3), we see that if $\vec{u}$ is continuous then $u$ satisfies the property of ultrahomogeneity with respect to any arrow between elements of the Fra\"iss\'e sequence $\vec{u}$; it remains to extend this result to hold for any arrow $f:a\rightarrow b$ in $\cal C$.              
So we have to prove that for any arrows $s:a\rightarrow u$ and $t:b\rightarrow u$ in $\cal D$ there exists an automorphism $\tilde{f}:u\rightarrow u$ such that $\tilde{f}\circ s=t\circ f$:
\[  
\xymatrix {
a \ar[d]_{f} \ar[r]^{s} & u \ar@{-->}[d]^{\tilde{f}}\\
b \ar[r]_{t} & u } 
\] 
Since $a$ is $\kappa$-small in $\cal D$, $s$ factors as $a\stackrel{s_{i}}{\rightarrow}u_{i}\stackrel{j_{i}}{\rightarrow}u$ for some $i\in \kappa$. From the fact that $u$ is $\cal C$-homogeneous (which we observed above), it follows that there exists an arrow $h:u_{i}\rightarrow u$ such that $h\circ s_{i}=t\circ f$. Now fact (3) implies the existence of an automorphism $\tilde{f}$ of $u$ such that $\tilde{f}\circ j_{i}=h$; then $\tilde{f}\circ s=\tilde{f}\circ j_{i}\circ {{s_{i}}}=h\circ s_{i}=t\circ f$, that is $\tilde{f}$ satisfies the required condition.\\
So far we have proved facts (1), (2) and (3).\\ 
Now, if $\cal C$ is $\kappa$-bounded, satisfies AP, JEP and has a dominating family of arrows $\cal F$ such that $|{\cal F}|\leq \kappa$, then by Theorem 3.5 in \cite{Kubis} there exists a $\kappa$-Fra\"iss\'e sequence in $\cal C$; hence by fact (2) there exists in ${\cal D}_{\kappa}$ a $\cal C$-homogeneous and $\cal C$-universal object $u$.\\
Let us now suppose that all the morphisms in ${\cal D}_{\kappa}^{c}$ are monic. If $u$ is a $\cal C$-universal and $\cal C$-homogeneous object in ${\cal D}_{\kappa}^{c}$ then by facts (1) and (2) above $u$ is $\cal C$-ultrahomogeneous. Now, suppose that $u,v\in {\cal D}_{\kappa}^{c}$ are both $\cal C$-ultrahomogeneous and $\cal C$-universal. Then by writing $u=\varinjlim \vec{u}$ and $v=\varinjlim \vec{v}$ where $\vec{u}$ and $\vec{v}$ are continuous inductive $\kappa$-sequences in $\cal C$, from fact (1) above we deduce that both $\vec{u}$ and $\vec{v}$ are continuous $\kappa$-Fra\"iss\'e sequences in $\cal C$; then by fact (3) there is an isomorphism $u\cong v$.\\
It remains to prove the last part of the theorem. From the proof of fact (2) above, it follows that there exists in $\cal C$ a $\kappa$-Fra\"iss\'e sequence satisfying the extension property; then the thesis follows from Proposition 3.1 in \cite{Kubis}.\\  
         
\end{proofs}
\begin{rmk}
\emph{One can relax the condition in the second part of the theorem that all the morphisms in ${\cal D}_{\kappa}^{c}$ are monic to the weaker condition that all the universal colimit arrows to the colimits of continuous $\kappa$-chains in $\cal D$ are monic in ${\cal D}_{\kappa}$, which is what one just needs in the proof of the theorem; however, in case all the morphisms in $\cal C$ are monic, this weaker condition turns out to be equivalent to the original condition.}
\end{rmk}
\begin{rmk}
\emph{The categorical theorem in Droste and G\"obel \cite{DG1} can be obtained as the particular case of Theorem \ref{teofond} when $i$ is the embedding of the category ${\cal C}_{\lt\kappa}$ of $\kappa$-small objects of a $\kappa$-algebroidal category $\cal C$ whose morphisms are all monic into $\cal C$ and $\cal F$ is the collection of arrows of some skeleton of ${\cal C}_{\lt\kappa}$. Fra\"iss\'e's theorem is already a particular case of Droste and G\"obel's result (as it is observed in \cite{DG}), hence it is \emph{a fortiori} a particular case of our theorem.}       
\end{rmk}
Let us note that given a category $\cal C$ as in Theorem \ref{teofond}, there is always an embedding $i:{\cal C}\rightarrow {\cal D}$ satisfying the hypotheses of the first part of the theorem, that is such that $\cal D$ has all colimits of $\kappa$-chains in $\cal C$ and all the objects in $\cal C$ are $\kappa$-small in $\cal D$; in fact, one can take as $\cal D$ the ind-completion $\Ind{\cal C}$ of $\cal C$ or the completion $(\Ind{\cal C})_{\kappa}$ of $\cal C$ in $\Ind{\cal C}$ under colimits of $\kappa$-chains. Recall that in case $\cal C$ is Cauchy-complete, $\cal C$ can be recoved from $\Ind{\cal C}$ as the full subcategory on the finitely presentable objects; also, we have seen above that $(\Ind{\cal C})_{\omega}$ can be identified with the full subcategory of $\Ind{\cal C}$ on the $\omega^{+}$-presentable objects.\\
Let us now apply Theorem \ref{teofond} in the context of first-order theories.
\begin{corollary}\label{cor_fo}
Let $\Sigma$ be a one-sorted signature, $T$ a first-order theory over $\Sigma$ and $\kappa$ an infinite cardinal such that $\kappa \gt card(\Sigma)$. Let ${{\mathbb T}\textrm{-mod}}_{e}$ be the category of $\mathbb T$-models and elementary embeddings between them and $i_{k}:{{\mathbb T}\textrm{-mod}}_{e}^{\kappa}\rightarrow  {{\mathbb T}\textrm{-mod}}_{e}$ be the embedding of the full subcategory ${{\mathbb T}\textrm{-mod}}_{e}^{\kappa}$ of ${{\mathbb T}\textrm{-mod}}_{e}$ on the $\kappa$-presentable objects into ${{\mathbb T}\textrm{-mod}}_{e}$. Then if ${{\mathbb T}\textrm{-mod}}_{e}^{\kappa}$ satisfies AP, JEP and has a dominating family of arrows in it of cardinality at most $\kappa$, $\mathbb T$ has a model of cardinality $\leq \kappa$ which is ${{\mathbb T}\textrm{-mod}}_{e}^{\kappa}$-ultrahomogeneous and ${{\mathbb T}\textrm{-mod}}_{e}^{\kappa}$-universal; moreover, a $\mathbb T$-model with these properties is unique (up to isomorphism) among the $\mathbb T$-models of cardinality $\leq \kappa$.    
\end{corollary}
\begin{proofs}
This immediately follows from Theorem \ref{teofond}, Proposition 1 in \cite{rosicky} and the remarks preceding Theorem \ref{teofond}.  
\end{proofs}
 
Finally, some cardinality considerations. If $\cal C$ is a category structured over $\Set$, or more generally over a functor category $[I,\Set]$ (where $I$ is a set, regarded here as a discrete category) via a ``forgetful'' functor $U:{\cal C}\rightarrow [I,\Set]$, then one can naturally define a notion of cardinality for objects of $\cal C$. Indeed, one can define the cardinality of an object $c\in \cal C$ by the formula $card(c)=|\coprod_{i\in I}U(c)(i)|=\coprod_{i\in I}|U(c)(i)|$. These definitions apply for instance to the case of models of a many-sorted (geometric) theory (in this case $\cal C$ is the category of such models while $I$ is the set of sorts of the theory), giving a notion of cardinality for such models that generalizes the definition of cardinality of a model in classical model theory. Suppose $i:{\cal C}\rightarrow{\cal D}$ is an embedding as in Theorem \ref{teofond}; if $\cal D$ is structured over a functor category $[I,\Set]$ via a functor $U:{\cal D}\rightarrow [I,\Set]$ then we have a notion of cardinality for objects of $\cal D$ and in particular of $\cal C$, and we might want to estimate the cardinality of the ultrahomogeneous universal object given by Theorem \ref{teofond} in terms of the cardinality of the objects of $\cal C$. This is particularly easy to do in case the functor $U$ creates colimits of $\kappa$-chains; in fact we know that the colimits in $[I,\Set]$ are computed pointwise and we have a particularly elegant description of filtered colimits (in particular colimits of $\kappa$-chains) in $\Set$ (see for example p. 77 in \cite{borceux}). Specifically, if $u=\varinjlim \vec{u}$ is the colimit in $\cal D$ of an inductive $\kappa$-sequence with values in $\cal C$, we have $card(u)=card(\varinjlim_{\cal D} \vec{u})=card(\varinjlim_{[I,Set]}(U\circ \vec{u}))=\coprod_{i\in I}|\varinjlim_{\Set}(U\circ \vec{u})(i)|$. Notice that for each $i\in I$, $(U\circ \vec{u})(i)$ defines a $\kappa$-chain in $\Set$. From this expression one can then deduce that if $|I|\leq \kappa$ and for each $i\in I$ and $j\in \kappa$ $|(U\circ \vec{u})(i)(j)|\leq \kappa$ then $card(u)\leq \kappa$. Thus for example if all the objects in $\cal C$ have cardinality $\leq \kappa$ and $|I|\leq \kappa$ then every object in ${\cal D}_{\kappa}$ has cardinality $\leq \kappa$. This is for instance the case of the classical Fra\"iss\'e's construction, where in fact the Fra\"iss\'e's limit is always at most countable.   
 
\section{The topos-theoretic interpretation}
A remark on notation: all the toposes in this section will be Grothendieck toposes, if not otherwise stated.\\
Let us recall that there exists an initial object in the category of toposes and geometric morphisms, which is given by the terminal category $1$ having just one object and the identity morphism on it; in fact, this category is a (coherent, atomic) Grothendieck topos, being the category of sheaves on the empty category with respect to the atomic topology on it (another presentation of it is obtained by taking the sheaves on $1$ with respect to the maximal Grothendieck topology on it, that is the topology in which all sieves cover). We will say that a topos $\cal E$ is trivial if it is naturally equivalent to $1$; of course, this is the same as saying that $\cal E$ is degenerate, that is $0_{\cal E}\cong 1_{\cal E}$.\\
Let us recall that a topos $\cal E$ is said to have enough points if the inverse image functors of the geometric morphisms $\Set\rightarrow \cal E$ are jointly conservative; every coherent topos has enough points (see for example \cite{El2}).
    
\begin{lemma}\label{teo1}
Let $\cal E$ be a topos with enough points. Then $\cal E$ is trivial if and only if it has no points.
\end{lemma}

\begin{proofs}
In one direction, let us suppose $\cal E$ trivial. Then $\cal E$ has no points because if $f:\Set\rightarrow \cal E$ were a point then we would have $0_{\Set}\cong f^{\ast}(0_{\cal E})\cong f^{\ast}(1_{\cal E})\cong 1_{\Set}$, which is absurd. Conversely, if $\cal E$ has no points then by taking the unique arrow $0:0_{\cal E}\rightarrow 1_{\cal E}$ in $\cal E$ then we trivially have that for each point $f$ of $\cal E$ $f^{\ast}(0)$ is an isomorphism; from the fact that $\cal E$ has enough points we can thus conclude that $0$ is an isomorphism, that is $\cal E$ is trivial.  
\end{proofs}

\begin{lemma}\label{teo2}
Let $\cal C$ be a category satisfying the right Ore condition, and $J_{at}$ the atomic topology on it. Then $\Sh({\cal C},J_{at})$ is trivial if and only if $\cal C$ is the empty category.
\end{lemma}

\begin{proofs}
Recall that $1_{\Sh({\cal C},J_{at})}$ is given by the constant functor $\Delta{1_{\Set}}:{\cal C}^{\textrm{op}}\rightarrow \Set$, while $0_{\Sh({\cal C},J_{at})}$ is given by the result of applying the associated sheaf functor $a:[{\cal C}^{\textrm{op}}, \Set]\rightarrow \Sh({\cal C},J_{at})$ to the initial object of $[{\cal C}^{\textrm{op}}, \Set]$, that is the constant functor $\Delta{\emptyset}:{\cal C}^{\textrm{op}}\rightarrow \Set$. But this functor is trivially a sheaf with respect to the atomic topology on $\cal C$, since all its covering sieves are non-empty, so $a(\Delta{\emptyset})\cong \Delta{\emptyset}$. Now, clearly, $\Delta{\emptyset}\cong \Delta{1_{\Set}}$ if and only if $\cal C$ is the empty category. 
\end{proofs}

\begin{lemma}\label{teo3}
Let $\cal C$ be a category satisfying the right Ore condition, and $J_{at}$ the atomic topology on it. Then if $[{\cal C}^{\textrm{op}}, \Set]$ is coherent, $\Sh({\cal C},J_{at})$ is coherent.
\end{lemma}

\begin{proofs}
From \cite{flatcoh} we know that if $[{\cal C}^{\textrm{op}}, \Set]$ is coherent, then we can axiomatize the theory of flat functors on $\cal C$ with coherent axioms in the language of presheaves on $\cal C$. Then, to obtain a coherent axiomatization for the theory of flat $J_{at}$-continuous functors on $\cal C$, it suffices to add to these axioms, for each arrow $f:c\rightarrow d$, the following (coherent) axiom:
\[
\top \: \vdash_{y}\: (\exists x\in c)(f(x)=y).
\]
\end{proofs}

{\flushleft
We recall that in \cite{flatcoh} Beke, Karazeris and Rosick\'y have introduced a notion of category having all fc finite limits and proved the following result: $[{\cal C}^{\textrm{op}}, \Set]$ is coherent if and only if $\cal C$ has all fc finite limits. Without going into details, we just remark that this fact can be profitably applied in connection with Lemma \ref{teo3} (see for example Theorem \ref{teocons} below).\\
 
We recall that a geometric theory $\mathbb T$ is said to be of presheaf type if its classifying topos is a presheaf topos (equivalently, the topos $[{\cal C},\Set]$, where ${\cal C}:=(\textrm{f.p.} {\mathbb T}\textrm{-mod}(\Set))$ is the category of finitely presentable $\mathbb T$-models in $\Set$). We will say that two geometric theories are Morita-equivalent if they have the same category of models - up to natural equivalence - into every Grothendieck topos $\cal{E}$ naturally in $\cal{E}$, equivalently the same classifying topos.\\
We recall from \cite{OC} that if $\mathbb T$ is a theory of presheaf type such that the category $(\textrm{f.p.} {\mathbb T}\textrm{-mod}(\Set))^{\textrm{op}}$ satisfies the right Ore condition (equivalently $\textrm{f.p.} {\mathbb T}\textrm{-mod}(\Set)$ satisfies AP), then the topos $\Sh((\textrm{f.p.} {\mathbb T}\textrm{-mod}(\Set))^{\textrm{op}}, J_{at})$ classifies the homogeneous $\mathbb T$-models. We note that the notion of homogeneity of a model of $\mathbb T$ in $\Set$ defined in \cite{OC} coincides with the notion of $(\textrm{f.p.} {\mathbb T}\textrm{-mod}(\Set))$-homogeneous object of the category ${\mathbb T}\textrm{-mod}(\Set)$ with respect to the embedding $(\textrm{f.p.} {\mathbb T}\textrm{-mod}(\Set))\hookrightarrow {\mathbb T}\textrm{-mod}(\Set)$ that we defined in the first section of this paper.\\ 
We will sometimes identify theories with their Morita-equivalence classes; the theory of flat $J_{at}$-continuous functors on $(\textrm{f.p.} {\mathbb T}\textrm{-mod}(\Set))^{\textrm{op}}$, which can be taken as the ``canonical'' representative for the Morita-equivalence class of theories classified by the topos $\Sh((\textrm{f.p.} {\mathbb T}\textrm{-mod}(\Set))^{\textrm{op}}, J_{at})$, will be called ``the theory of homogeneous $\mathbb T$-models''.}\\
A geometric theory is said to be consistent if it has at least one model in $\Set$.\\  
%{\flushleft 
The previous lemmas combine to give the following consistency result.%}

\begin{theorem}\label{teocons}
Let $\mathbb T$ be a theory of presheaf type such that the category $\textrm{f.p.} {\mathbb T}\textrm{-mod}(\Set)$ has the amalgamation property. If the theory of homogeneous $\mathbb T$-models is Morita-equivalent to a coherent theory (for example when the category $\textrm{f.p.} {\mathbb T}\textrm{-mod}(\Set)$ has all fc finite colimits) and there is at least one $\mathbb T$-model in $\Set$, then there exists at least one homogeneous $\mathbb T$-model in $\Set$.  
\end{theorem}

\begin{proofs}
The theory $\mathbb T'$ of homogeneous $\mathbb T$-models is Morita-equivalent to a coherent theory if and only if its classifying topos $\Sh((\textrm{f.p.} {\mathbb T}\textrm{-mod}(\Set))^{\textrm{op}}, J_{at})$ is a coherent topos. Notice that for any category $\cal C$, $\cal C$ is empty if and only if $\Ind{{\cal C}}$ is empty; so if $\mathbb T$ is a theory of presheaf type then $\mathbb T$ has a model in $\Set$ if and only if it has a finitely presentable model in $\Set$. Then, since $\textrm{f.p.} {\mathbb T}\textrm{-mod}(\Set)$ is not the empty category, it follows from Lemma \ref{teo2} that the topos $\Sh((\textrm{f.p.} {\mathbb T}\textrm{-mod}(\Set))^{\textrm{op}}, J_{at})$ is not trivial. Hence, by Lemma \ref{teo1}, it has a point. This point corresponds to a $\mathbb T'$-model in $\Set$, that is, to a homogeneous $\mathbb T$-model in $\Set$.
The fact that when the category $\textrm{f.p.} {\mathbb T}\textrm{-mod}(\Set)$ has all fc finite colimits, $\mathbb T'$ is Morita-equivalent to a coherent theory follows from Lemma \ref{teo3}.          
\end{proofs} 
A (many-sorted) geometric theory is said to be atomic if it is classified by an atomic topos. Of course, the property of atomicity for a theory is stable under Morita-equivalence. A geometric theory $\mathbb T$ over a signature $\Sigma$ is said to be complete if every sentence over $\Sigma$ is $\mathbb T$-provably equivalent to $\top$ or $\bot$, but not both.  It is well-known that if $\mathbb T$ is atomic then $\mathbb T$ is complete if and only if its classifying topos $\Set[\mathbb T]$ is connected (equivalently, two-valued - see the proof of Theorem \ref{teo4} below). Recall that if a theory is coherent then its completeness implies its consistency (cfr. for example Lemma \ref{teo1}), but this implication does not hold for a general geometric theory; in fact, there exist connected atomic toposes without points (see for example \cite{El2}). We also remark that the property of completeness for a geometric theory is stable under Morita-equivalence, being equivalent to a categorical property (to be two-valued) of the corresponding classifying topos.\\
 
\begin{theorem}\label{teo4}
Let $\cal C$ be a category satisfying the right Ore condition, and $J_{at}$ the atomic topology on it. Then the atomic topos $\Sh({\cal C},J_{at})$ is connected if and only if $\cal C$ is a connected category. 
\end{theorem}

\begin{proofs}
Recall that a topos $\cal E$ is said to be locally connected if the geometric morphism $\gamma:{\cal E}\rightarrow \Set$ is essential, that is the inverse image functor $\gamma^{\ast}:\Set \rightarrow \cal E$ has a left adjoint $\gamma_{!}:{\cal E}\rightarrow \Set$. An object $A$ of a locally connected topos $\cal E$ is said to be connected if $\gamma_{!}(A)\cong 1_{\Set}$. Every atomic topos $\cal E$ is locally connected (see for example p. 684 of \cite{El2}), and the objects of $\cal E$ which are connected are also called atoms.\\
We observe that an object $A$ of an atomic topos $\cal E$ is an atom if and only if the only subobjects of $A$ in $\cal E$ are $0_{A}:0\rightarrow A$ and $1_{A}:A\rightarrow A$ and they are distinct from each other. Indeed, this easily follows from the bijection $\Sub_{\cal E}(A)\cong \Sub_{\Set}(\gamma_{!}(A))$ (cfr. p. 685 of \cite{El2}). Hence, since every atomic topos is locally connected, Lemma C.3.3.3 in \cite{El2} gives the following characterization, to which we refer as to $(\ast)$: an atomic topos $\cal E$ is connected if and only if the only subobjects of $1_{\cal E}$ in $\cal E$ are $0_{1}:0\rightarrow 1$ and $1_{1}:1\rightarrow 1$ and they are distinct from each other. We use this criterion to prove our theorem.\\ 
We can identify the subterminals in $\Sh({\cal C},J_{at})$ with $J_{at}$-ideals on $\cal C$ (see p. 576 of \cite{El2}). By recalling (from the proof of Lemma \ref{teo2}) that $0_{\Sh({\cal C},J_{at})}$ is the constant functor $\Delta{\emptyset}:{\cal C}^{\textrm{op}}\rightarrow \Set$, condition $(\ast)$ can thus be rephrased as follows:\\
$\cal C$ is non-empty and every non-empty subset $I\subseteq ob({\cal C})$ which is a sieve (that is, for each arrow $f:a\rightarrow b$ in $\cal C$, $b\in I$ implies $a\in I$) and satisfies the property $(\forall R\in J_{at}(U)((\forall f_{i}:U_{i}\rightarrow U \in R, U_{i}\in I)\imp (U\in I))$ is the whole of $\cal C$.\\
Being $J_{at}$ the atomic topology on $\cal C$, this condition simplifies to:\\
$\cal C$ is non-empty and every non-empty subset $I\subseteq ob({\cal C})$ which is a sieve and satisfies the property $\forall f:V\rightarrow U$ in $\cal C$, $((V\in I) \imp (U\in I))$ is the whole of $ob(\cal C)$; but this is clearly equivalent to saying that $\cal C$ is connected.  
\end{proofs}

\begin{theorem}\label{teo5}
Let $\cal C$ be a non-empty category satisfying the amalgamation property. Then $\cal C$ satisfies the joint embedding property if and only if it is a connected category. 
\end{theorem}

\begin{proofs}
If $\cal C$ satisfies JEP then for any objects $a,b\in {\cal C}$ there exists an object $c\in \cal C$ and morphisms $f:a\rightarrow c$, $g:b\rightarrow c$ in $\cal C$:
\[  
\xymatrix {
 & a  \ar[d]^{f} \\
b \ar[r]_{g} & c } 
\] 
Then we have the following zig-zag between $a$ and $b$:
\[  
\xymatrix {
 & a \ar[dl]_{1_{a}} \ar[dr]^{f} &  & b \ar[dl]_{g} \ar[dr]^{1_{b}} \\
a & & c & & b. } 
\] 
Conversely, we prove that for any objects $a,b\in {\cal C}$ there exists an object $c\in \cal C$ and morphisms $f:a\rightarrow c$, $g:b\rightarrow c$ in $\cal C$ by induction on the length $n$ of a zig-zag that connects $a$ and $b$. If $n=1$ then the thesis follows immediately from the amalgamation property. If $n\gt 1$ we have a zig-zag 
\[  
\xymatrix {
 & \ldots & \ldots  & d'_{n} \ar[dl]_{f_{n}} \ar[dr]^{g_{n}} \\
d_{0}=a & \ldots & d_{n-1} & & d_{n}=b. } 
\] 
By applying the induction hypothesis to the pair $a, d_{n-1}$ one gets an object $d\in \cal C$ and morphisms $h:a\rightarrow d$, $k:d_{n-1}\rightarrow d$ in $\cal C$. The amalgamation property applied to the pair of morphisms $k\circ f_{n}$ and $g_{n}$ then gives an object $c$ and two morphisms $s:d\rightarrow c$ and $t:b\rightarrow c$. Then we have morphisms $f:=s\circ h:a\rightarrow c$ and $g:=t:b\rightarrow c$, as required.
\end{proofs}

From Theorems \ref{teo4} and \ref{teo5} we thus deduce that given a consistent theory of presheaf type $\mathbb T$ such that the category $\textrm{f.p.} {\mathbb T}\textrm{-mod}(\Set)$ satisfies the amalgamation property, the condition that $\textrm{f.p.} {\mathbb T}\textrm{-mod}(\Set)$ satisfies JEP is exactly what makes the theory $\mathbb T'$ of homogeneous $\mathbb T$-models complete. Indeed, $\mathbb T'$ is complete if and only if $\Set[\mathbb T']\simeq \Sh((\textrm{f.p.} {\mathbb T}\textrm{-mod}(\Set))^{\textrm{op}}, J_{at})$ is connected, if and only if $(\textrm{f.p.} {\mathbb T}\textrm{-mod}(\Set))^{\textrm{op}}$ is connected, if and only if $(\textrm{f.p.} {\mathbb T}\textrm{-mod}(\Set))$ is connected, if and only if $(\textrm{f.p.} {\mathbb T}\textrm{-mod}(\Set))$ satisfies JEP.\\

A geometric theory is said to be countably categorical if any two countable models of $\mathbb T$ in $\Set$ are isomorphic (where by `countable' we mean either finite or denumerable). Notice that, by our definition, any geometric theory having no models in $\Set$ is (vacously) countably categorical. We recall from \cite{OC5} that every complete atomic geometric theory is countably categorical; so, by the remarks above, we obtain the following result.

\begin{theorem}\label{teo6}
Let $\mathbb T$ be a consistent theory of presheaf type such that the category $\textrm{f.p.} {\mathbb T}\textrm{-mod}(\Set)$ has the amalgamation and joint embedding properties. If $\mathbb T'$ is a geometric theory which is Morita-equivalent to the theory of homogeneous $\mathbb T$-models then $\mathbb T'$ is complete and countably categorical.
\end{theorem}
\qed

\begin{rmk}
\emph{Concerning the existence of homogeneous $\mathbb T$-models in $\Set$, we note that if the theory $\mathbb T$ in Theorem \ref{teo6} is coherent then, by Lemma \ref{teo3} and Theorem \ref{teocons}, there exists a homogeneous $\mathbb T$-model in $\Set$. If moreover the signature of $\mathbb T$ is countable then, by the results in \cite{OC5}, there is a countable homogeneuos $\mathbb T$-model in $\Set$.}  
\end{rmk}

The usefulness of Theorem \ref{teo6} lies in the fact that it is generally not difficult to see, given a theory of presheaf type $\mathbb T$, if a certain theory is Morita-equivalent to the theory of homogeneous $\mathbb T$-models. In fact, one can use Corollary 4.7 in \cite{OC} and the explicit description of the homogeneous models given in \cite{OC}. For example, in \cite{OC} we saw that, given the theory $\mathbb T$ of linear orders, the dense linear orders without endpoints corresponded precisely to the homogeneous $\mathbb T$-models. By using similar methods, one can also show that, given the theory of decidable objects, the infinite decidable objects are exactly the homogeneous decidable objects and that, given the algebraic theory of Boolean algebras, the atomless Boolean algebras are exactly the homogeneous Boolean algebras. This leads, via Theorem \ref{teo6}, to an alternative proof that the theory of dense linear orders without endpoints and the theory of atomless Boolean algebras are complete and countably categorical.\\  
Moreover, we know from \cite{OC5} that, under the hypotheses of Theorem \ref{teo6}, the Booleanization of the theory $\mathbb T$ axiomatizes the $\mathbb T$-homogeneous models, and hence we may deduce that any two countable $\mathbb T$-homogeneous models in $\Set$ are isomorphic (cfr. Theorem 3.3 \cite{OC5}).\\  
We also recall from \cite{OC5} that if $\mathbb T$ is an atomic, complete countable geometric theory with infinite models in $\Set$ then, denoted by $M$ the unique countable model of $\mathbb T$ (up to isomorphism), we have the following representation for the classifying topos $\Set[\mathbb T]$ of $\mathbb T$:
\[
\Set[\mathbb T]\simeq \Cont(Aut(M)),
\]
where $\Cont(Aut(M))$ is the topos of continuous $Aut(M)$-sets, $Aut(M)$ being endowed with the topology of pointwise convergence.\\ 
      
Let us now apply the categorical theorem in the first section in the context of the theories of presheaf type.
\begin{theorem}
Let $\mathbb T$ be a consistent theory of presheaf type such that the category $\textrm{f.p.} {\mathbb T}\textrm{-mod}(\Set)$ satisfies AP, JEP. If there exists in $\textrm{f.p.} {\mathbb T}\textrm{-mod}(\Set)$ a dominating family of arrows of finite or countable cardinality then there exists in $\Set$ a $\textrm{f.p.} {\mathbb T}\textrm{-mod}(\Set)$-homogeneous and $(\textrm{f.p.} {\mathbb T}\textrm{-mod}(\Set))$-universal ${\mathbb T}$-model; also, given a $(\textrm{f.p.} {\mathbb T}\textrm{-mod}(\Set))$-homogeneous and $(\textrm{f.p.} {\mathbb T}\textrm{-mod}(\Set))$-universal ${\mathbb T}$-model $M$, if $M$ can be written as the colimit in ${\mathbb T}\textrm{-mod}(\Set)$ of a $\omega$-chain of finitely presentable $\mathbb T$-models (equivalently, is $\omega^{+}$-presentable) then, provided that all the morphisms in $({\mathbb T}\textrm{-mod}(\Set))_{\omega}$ are monic, $M$ is $(\textrm{f.p.} {\mathbb T}\textrm{-mod}(\Set))$-ultrahomogeneous and unique (up to isomorphism) with this property among the $(\textrm{f.p.} {\mathbb T}\textrm{-mod}(\Set))$-universal and $\omega^{+}$-presentable $\mathbb T$-models.\\ If $\mathbb T'$ is a geometric theory whose models (in any Grothendieck topos) are the homogeneous $\mathbb T$-models, then $\mathbb T'$ is complete and countably categorical. In particular, if $\mathbb T$ is countable and has infinite models in $\Set$ there exists a unique (up to isomorphism) countable homogeneous $\mathbb T$-model $M$, and 
\[
\Set[\mathbb T']\simeq \Sh((\textrm{f.p.} {\mathbb T}\textrm{-mod}(\Set))^{\textrm{op}}, J_{at}) \simeq \Cont(Aut(M)),
\]
$Aut(M)$ being endowed with the topology of pointwise convergence.  
\end{theorem}

\begin{proofs}
This is immediate from Theorem \ref{teofond}, the remarks following it, Theorem \ref{teo6} and the remark above.
\end{proofs}

Let us now introduce the following notions.\\
Given an embedding $i:{\cal C}\hookrightarrow {\cal D}$ and an object $u\in {\cal D}$ together with a choice of an arrow $f_{c}:c\rightarrow u$ in $\cal D$ for each object $c$ of $\cal C$, we can consider a category $\tilde{\cal C}$, defined as the full subcategory of $({\cal C}\downarrow u)$ on the arrows $f:a\rightarrow u$ in $\cal D$ such that there exists an automorphism $\alpha$ of $u$ (that is, an isomorphism $\alpha:u\rightarrow u$ in the category $\cal D$) such that $f=\alpha \circ f_{a}$. Then we can define a functor $\chi:\tilde{\cal C}^{\textrm{op}}\rightarrow Subgr(Aut(u))$, where $Subgr(Aut(u))$ is the collection of the subgroups of $Aut(u)$ regarded as a poset category with respect to the inclusion, in the following way: $\chi$ sends an object $f:a\rightarrow u$ in $\tilde{\cal C}$ to the subgroup $Aut_{f}(u)$ of $Aut(u)$ formed by the automorphisms $\alpha$ of $u$ such that $\alpha\circ f=f$ and an arrow $h:f\rightarrow g$ in $\tilde{\cal C}$ to the inclusion $Aut_{g}(u)\subseteq Aut_{f}(u)$. If $\chi$ is full and faithful and reflects identities (that is, for each pair of arrows $h,k$ in $\tilde{\cal C}^{\textrm{op}}$ $\chi(h)=\chi(k)$ implies $h=k$) we say that $u$ satisfies the Galois property with respect to $\tilde{\cal C}$; notice that if $\cal C$ is skeletal and $\chi$ is full and faithful then $u$ satisfies the Galois property with respect to $\tilde{\cal C}$. Also, we can endow the group $Aut(u)$ with a topology $\cal U$ by saying that the subgroups in the image of the functor $\chi$ form a base of neighbourhoods of the identity.\\     
In the context of these notions, the following proposition holds.
\begin{proposition}
Given an embedding $i:{\cal C}\hookrightarrow {\cal D}$ such that all the arrows $f_{c}$ (for $c\in {\cal C}$) are monic, let $u$ be a $\cal C$-ultrahomogeneous object in $\cal D$ which satisfies the Galois property with respect to $\tilde{\cal C}$. Then the category $\cal C$ satisfies the amalgamation property and there is a natural equivalence
\[
\Sh({\cal C}^{\textrm{op}},J_{at})\simeq \Cont(Aut(u))
\] 
where $\Cont(Aut(u))$ is the topos of continuous $Aut(u)$-sets, $Aut(u)$ being endowed with the topology $\cal U$.
\end{proposition}
\begin{proofs}
From Theorem 2 p. 154 \cite{MM} we deduce that $\Cont(Aut(u))$ is naturally equivalent to $\Sh({\bf S}_{\cal U}(Aut(u)),J_{at})$, where ${\bf S}_{\cal U}(Aut(u))$ is the category having as objects the continuous $Aut(u)$-sets of the form $Aut(u)\slash \chi(f)$ for $f\in \tilde{\cal C}$ and as arrows $Aut(u)\slash \chi(f)\rightarrow Aut(u)\slash \chi(g)$ the cosets $\chi(g)\alpha$ with the property that $\chi(f)\subseteq \alpha^{-1}\chi(g)\alpha$ (see \cite{MM} for more details). To prove our proposition it is therefore enough to show that there is an equivalence of categories between ${\bf S}_{\cal U}(Aut(u))$ and ${\cal C}^{\textrm{op}}$. We explicitly define a functor $F:{\bf S}_{\cal U}(Aut(u))\rightarrow {\cal C}^{\textrm{op}}$ and prove that it is an equivalence of categories.
Let us first define $F$ on objects: $F$ sends an object $Aut(u)\slash \chi(f)$ of ${\bf S}_{\cal U}(Aut(u))$ to $dom(f)\in {\cal C}$; this is well-defined since $\chi$ reflects identities. Given an arrow $Aut(u)\slash \chi(f)\rightarrow Aut(u)\slash \chi(g)$, represented by a coset $ \chi(g)\alpha$, we have that $\chi(f)\subseteq \alpha^{-1}\chi(g)\alpha$, equivalently $\alpha \chi(f)\alpha^{-1}\subseteq \chi(g)$. This means that $\alpha\circ \beta \circ \alpha^{-1}\circ g = g$ (equivalently, $\beta \circ (\alpha^{-1}\circ g)=(\alpha^{-1}\circ g)$) for each $\beta\in Aut(u)$ such that $\beta \circ f=f$, which is in turn equivalent to saying that $\chi(f)\subseteq \chi(\alpha^{-1}\circ g)$. This implies, by our hypothesis that $\chi$ is full and faithful, that there exists a unique arrow $z:dom(g)\rightarrow dom(f)$ in $\cal C$ such that $f\circ z=\alpha^{-1}\circ g$. We put $F(\chi(g)\alpha)=z$; this is well-defined since $\chi(g) \alpha =\chi(g) \alpha'$ if and only if $\alpha \circ \alpha'^{-1}\in \chi(g)$, if and only if $\alpha^{-1} \circ g=\alpha'^{-1} \circ g$, if and only if $f\circ F(\chi(g) \alpha)=f\circ F(\chi(g)\alpha')$ if and only if $F(\chi(g)\alpha)=F(\chi(g)\alpha')$, where the last equivalence follows from the fact that $f$ is monic. This also proves that $F$ is faithful. $F$ is full because $u$ is ${\cal C}$-ultrahomogeneous, and it is surjective by definition of $\cal U$. Therefore, $F$ is an equivalence of categories.          
\end{proofs}

\vspace{7 mm}
{\bf Acknowledgements.} I am very grateful to my Ph.D. supervisor Peter Johnstone for his support and encouragement. Thanks also to Martin Hyland for having suggested me to investigate Fra\"iss\'e's construction topos-theoretically.\\ 

\end{document}